# $L^P$-VARIATIONS FOR MULTIFRACTAL FRACTIONAL RANDOM WALKS

By Carenne Ludeña

*Instituto Venezolano de Investigaciones Cientificas*

A multifractal random walk (MRW) is defined by a Brownian motion subordinated by a class of continuous multifractal random measures $M[0,t]$, $0 \leq t \leq 1$. In this paper we obtain an extension of this process, referred to as multifractal fractional random walk (MFRW), by considering the limit in distribution of a sequence of conditionally Gaussian processes. These conditional processes are defined as integrals with respect to fractional Brownian motion and convergence is seen to hold under certain conditions relating the self-similarity (Hurst) exponent of the fBm to the parameters defining the multifractal random measure $M$. As a result, a larger class of models is obtained, whose fine scale (scaling) structure is then analyzed in terms of the empirical structure functions. Implications for the analysis and inference of multifractal exponents from data, namely, confidence intervals, are also provided.

**1. Introduction.** Multifractal models, or more generally, multifractal random walks, with respect to ordinary Brownian motion or fractional Brownian motion have become very popular models for scaling phenomena. This kind of phenomena is associated to very nonlinear models as encountered, for example, in hydrodynamic turbulence, biological rhythms or the network traffic [4, 5, 20, 21].

The scaling property is usually defined in terms of a scaling exponent, $\zeta(p)$, associated to the process $X(t)$ describing the phenomena. This exponent is defined as

$$\zeta(p) := \frac{\log \mathbb{E}(|X(t+\tau) - X(t)|^p)}{\log(\tau)}$$

at different scales $0 < \tau < T$. Based on the multifractal formalism [2, 14, 19] which establishes the relationship between $\zeta(p)$ and the characterization of the fluctuations of the local regularity of $X(t)$, estimating this scaling function has attracted a great deal of attention both from a theoretical and a practical point of view (see, e.g., [1, 16, 23]).

Typically, given a sample path of $X(t)$, $t \in [0,1]$, say, estimators of $\zeta(p)$ are based on the logarithm of the empirical structure function of the process over a certain regular partition $t_0, \ldots, t_m$ of $[0,1]$ ($t_i - t_{i-1} = \tau$), that is, of

$$S_m(p) := \frac{1}{m} \sum_{k=1}^{m} |X(t_k) - X(t_{k-1})|^p.$$

This procedure implicitly assumes that $S_m(p)$ converges a.s. to $\mathbb{E}(|X(t+\tau) - X(t)|^p)$ as $m \to \infty$ (i.e., ergodicity), and thus, that $\log(S_m(p))/\log(\tau)$ is a consistent estimator of $\zeta(p)$. We will refer to this method as a moment or $L^p$-variation estimator.

$L^p$-variation estimators have two major drawbacks in practice. On the one hand, they tend to be very data demanding [23]. On the other hand, out of a certain range $[p_-, p_+]$ a certain linearization effect ([16, 22, 23] and the references therein) tends to appear. That is, $\log S_m(p)/\log(\tau) \to kp$ for a certain constant $k$ instead of the expected $\zeta(p)$. This has led to consider estimators based not on the increments, but rather on the suprema of the increments or likewise techniques which provide more robust methods [2, 16, 23], although this phenomena is not thoroughly understood.

The aim of this paper is twofold. In the first place, based on the construction of multifractal random measures (MRM) developed by Bacry and Muzy in [7], we consider an extension of their multifractal random walks for ordinary Brownian motion to multifractal fractional random walks (MFRW) based on fractional Brownian motion with Hurst parameter $H$. Two possible courses of action are proposed. In each case the existence of the process is assured under certain additional conditions involving the Hurst parameter and the multifractal random measure.

In the second place, we study the conditional distribution, given the multifractal random measure, of the structure function for even and positive $p$. As a corollary, we discuss the properties of a $L^p$-variation estimator of the scaling function, for even $p$. The stated conditional asymptotic distribution allows us to construct confidence intervals for $\zeta(p)$ for large values of $m$, the number of observations. Our main result shows that approximation errors are fundamentally due to a bias term which tends to zero at a logarithmic rate.

On the other hand, a very interesting phenomena, related to a certain linearization effect reported in many applications, appears for intermediate values of $m$. This phenomena is related to the range of scales over which the scale invariance property holds true. In the MRM model developed in [7] this range depends on a certain parameter $T$. Based on this characterization, we show that the behavior of the estimator will be very different according to the relative size of $T$ and $m$. If $T \ll m$ ($T = 1$, say), the asymptotic behavior of the estimator will depend on the conditional asymptotic distribution of the structure function. However, if $T \sim m$, for big enough $m$, the scaling property will introduce a linear bias term in the estimation procedure.

The article is organized as follows. In Section 2 we introduce the class of MRM following the construction of Bacry and Muzy and define the MFRW based on fractional Brownian motion with Hurst parameter $H > 1/2$. In Section 3 we study the asymptotic properties of the sums of certain functions of the increments of the MFRW for $1/2 < H < 3/4$. In Section 4 we apply the asymptotic results of Section 3 to study the properties of a moment based estimator of the scaling function. In Section 4.2 we discuss the problem of estimating $\zeta(p)$ when $T \sim m$, that is, when the number of observations is comparable to the scale invariance range.

**2. Multifractal random measures.** The existence and analysis of discrete multifractal random measures, based on the Mandelbrot multiplicative cascades [17], was first established by Kahane and Peyriere [15]. Continuous extensions were introduced by Barral and Mandelbrot [8]. This class of models was further extended by Bacry and Muzy [7] who define a class of stationary random measures (MRM) based on the following elements:

1. Let $\nu$ be a *Lévy measure* defined over $\mathbb{R}$ satisfying

$$\left(\int_{-\infty}^{-u} + \int_{u}^{\infty}\right) \frac{\nu(dx)}{x^2} < \infty,$$

for all $u > 0$. Associated to $\nu$ and $m \in \mathbb{R}$, define for $q \in \mathbb{R}$

$$\phi(q) = imq + \int \frac{e^{iqx} - 1 - iq\sin x}{x^2} \nu(dx).$$

2. Set $S^+ = \{(t, l), t \in \mathbb{R}, l \in \mathbb{R}^+\}$ and the associated measure

$$\mu(dt, dl) = l^{-2} \, dt \, dl.$$

3. Define the process $P(A), A \subset S^+$, such that, for any sequence of disjoint sets $A_n \subset S^+$, $(P(A_n))_{n=1}^{\infty}$ are independent r.v., $P(\bigcup_{n=1}^{\infty} A_n) = \sum_{n=1}^{\infty} P(A_n)$ a.s. and

$$\mathbb{E}(e^{iqP(A)}) = e^{\phi(q)\mu(A)}.$$

4. Let

$$f(k) = \begin{cases} k, & k \leq T, \\ T, & k > T, \end{cases}$$

where $T$ is arbitrarily large. For $t \in \mathbb{R}$ and $l > 0$, define the collection of sets

(1) $$A_l(t) = \{(s,k), k \geq l, -f(k)/2 < s - t < f(k)/2\}.$$

Next for $t \in \mathbb{R}$, consider the process $w_l(t) = P(A_l(t))$.

5. For $l > 0$ and $I$ a Lebesgue measurable set, define

$$M_l(I) = \int_I e^{w_l(t)} \, dt,$$

and set $M_l(dt)$ to be the associated measure. It can be seen (see exact conditions and references below) that there exists a random measure $M(dt)$ such that a.e. $M_l(dt)$ converges weakly toward $M(dt)$. The MRM measure is then defined as the limit

$$M(dt) = \lim_{l \to 0^+} M_l(dt).$$

The properties, and existence, of $M(dt)$ are discussed in [7], as well as several examples of popular choices of the measure $\nu$.

Associated to the set process $P$ and $\phi$, we may also consider the real function $\psi(q) = \phi(-iq)$ (see [7] for further details) defined in terms of the moment generating function, $\mathbb{E}(e^{qP(A)}) = e^{\psi(q)\mu(A)}$, for any $A \subset S^+$. Throughout the article we shall assume $\psi(1) = 0$. Let $\zeta(q) = q - \psi(q)$, and introduce the following assumptions:

[A1] There exists $\varepsilon > 0$ such that $\zeta(1 + \varepsilon) > 1$
[A(q)] $\zeta(q) > 1$.

The following lemma provides a summary of the properties of MRM obtained in [7] and which we list here for ease of reference.

LEMMA 2.1 (Theorems 2, 3, 4 and Lemma 5 in [7]).

1. *[Nondegeneracy of $M(dt)$.]* If there exists $\varepsilon > 0$ such that $\zeta(1+\varepsilon) > 1$ ([A1]), then $\mathbb{E}(M([0,t])) = t$.
2. *(Positive order moments of $M$.)* If $\zeta(q) > 1$ ([A(q)]) for $q > 0$, then $\mathbb{E}(M([0,t])^q) < \infty$. Moreover, if [A1] holds and $\mathbb{E}(M([0,t])^q) < \infty$, then $\zeta(q) \geq 1$.
3. *[Moments of $\sup_{u \in [0,t]} e^{w_l(u)}$.]* If $\psi(q) \neq \infty$ [in particular, if $\zeta(q) > 1$], then

$$\mathbb{E}\left(\sup_{u \in [0,t]} e^{qw_l(u)}\right) < \infty.$$

4. *(Scaling property.)* For all $t \leq T$,
$$\mathbb{E}(M[0,t]^q) = \left(\frac{t}{T}\right)^{\zeta(q)} \mathbb{E}(M[0,T]^q),$$
whenever $\zeta(q) > 1$. Moreover,
$$\{w_{l\lambda}(\lambda t)\}_t \stackrel{\text{law}}{=} \Omega_\lambda + \{w_l(t)\}_t$$
and
$$\{M[0, \lambda t]\}_t \stackrel{\text{law}}{=} \lambda e^{\Omega_\lambda} \{M[0,t]\}_t$$
for all $\lambda \in (0,1)$, where $\Omega_\lambda$ is an infinitely divisible random variable which is independent of $\{w_l(t)\}_{l,t}$ with characteristic function
$$\mathbb{E}(e^{iq\Omega_\lambda}) = \lambda^{-\phi(q)}.$$

2.1. *Constructing multifractal fractional random walks.* Assume throughout this subsection, without loss of generality and in order to simplify the construction of the MFRW, that $T = 1$. For $M$ a nondegenerate measure (i.e., such that $\mathbb{E}(M([0,t])) = t$), Bacry and Muzy [7] define the subordinate random walk as
$$Y_t = B(M([0,t])),$$
for $B(t)$ a Brownian motion which is independent of $M$. This process can be interpreted as the centered, conditionally Gaussian process whose conditional covariance function is defined as the limit
$$\Gamma(s,t) := \lim_{l \to 0^+} \int_0^{t \wedge s} e^{w_l(u)} \, du.$$

The authors also show this process is equivalent in law to the limit process
$$Y'_t = \lim_{l \to 0^+} \int_0^t e^{w_l(s)/2} \, dB_s.$$

This procedure has been repeated based on heuristic arguments for fractional Brownian motion (see, e.g., [6, 16]). The work by Bacry and Muzy [7] provides the basis for a rigorous definition based on the properties of the stochastic integral with respect to fBm whenever the integrand is deterministic (see, e.g., [10]) and the Hurst coefficient, $H$, is greater than $1/2$. However, though it might appear to be the natural thing to do, the MFRW cannot be defined directly as the limit in distribution of the (well-defined) stochastic integrals $\int_0^t e^{w_l(s)/2} \, dB_s^H$. Indeed, although for each $l > 0$, the conditional covariance $R_l(s,t)$ associated to the stochastic integral, with respect to fBm for $H > 1/2$, is well defined and given by
$$R_l(s,t) = C_H \int_0^t \int_0^s |u-v|^{2H-2} e^{w_l(u)/2} e^{w_l(v)/2} \, du \, dv,$$

a direct calculation of the expectation using the moment generating function of process $\{w_l\}_t$ (Lemma 1 in [7]), yields that $\liminf_l \mathbb{E} R_l(s,t) \to 0$. More precisely (assume $s < t$), we have

$$\mathbb{E} R_l(s,t) = e^{2\psi(1/2)} l^{-2\psi(1/2)} \int_0^t \int_0^s |u-v|^{2H-2} e^{-2\psi(1/2)\rho_l(u-v)} \, du \, dv$$

$$= \left( 2 \int_0^s \int_0^u v^{2H-2} e^{-2\psi(1/2)\rho_l(v)} \, dv \, du \right.$$

$$\left. + \int_s^t \int_{u-s}^u v^{2H-2} e^{-2\psi(1/2)\rho_l(v)} \, dv \, du \right) / l^{2\psi(1/2)}$$

$$= l^{-2\psi(1/2)} 2 \left[ \int_0^l \int_0^u v^{2H-2} e^{-2\psi(1/2)\rho_l(v)} \, dv \, du \right.$$

$$\left. + \int_l^s \left[ \int_0^l + \int_l^u \right] v^{2H-2} e^{-2\psi(1/2)\rho_l(v)} \, dv \, du \right]$$

$$+ e^{2\psi(1/2)} l^{-2\psi(1/2)} \int_s^t \int_{u-s}^u v^{2H-2} e^{-2\psi(1/2)\rho_l(v)} \, dv \, du$$

$$= I + II + III + IV,$$

where

$$\rho_l(v) = \begin{cases} \ln(1/l) + 1 - v/l, & \text{if } v \leq l, \\ \ln(1/v), & \text{if } l \leq v \leq 1, \\ 0, & \text{if } v > 1. \end{cases}$$

Setting $a = 2\psi(1/2)$

$$I = 2l^{2H} \left( \int_0^l t^{2H-2}(t-1) e^{at} \, dt \right) \to 0$$

and

$$II = 2(s-l) l^{2H-1} \int_0^1 t^{2H-2}(t-1) e^{at} \, dt \to 0$$

On the other hand, for $l < v < 1$, we have $\rho_l(v) = \log(1/v)$, so that

$$III \leq 2e^a \left( l^{-2\psi(1/2)} \int_0^s t^{2H-2+2\psi(1/2)+1} \, dt - l^{2H-1} \right) \to 0$$

since $H > 1/2$ and $\psi(1/2) \leq 0$ by convexity. Finally, $IV$ can be bounded as above to show it also tends to zero.

Set $R(s,t) := \lim_{l \to 0^+} R_l(s,t)$. An application of Fatou's lemma yields $\mathbb{E} R(s,t) = 0$, and whence that $R(s,t) \equiv 0$.

In what follows we propose two alternative procedures. In each case an additional condition relating $H$ and the scaling function $\psi(\cdot)$ must be imposed in order to assure the existence of the limiting process.

2.1.1. *A first approach.* Start by defining the measure MRM $M^{1/2}$, whenever the scaling function $\psi(1/2) > -1/2$, as follows. For any Lebesgue measurable set $I$, $M^{1/2}(I) = \lim_{l \to 0^-} l^{-\psi(1/2)} e^{-\psi(t)/2} e^{w_l(t)/2} \, dt$. Existence of this measure follows as for the proof of Theorem 1 in [7].

With this construction, we have the following definition of a MFRW, where here and below $\mathbb{E}_M$ denotes the conditional expectation given the MRM $M$.

DEFINITION 2.1. Let $H > 1/2$. The MFRW-1/2 is the centered process $X_t$ which, conditionally to $M(dt)$, has a Gaussian distribution with conditional covariance

$$S(s,t) = C_H \int_0^t \int_0^s |u-v|^{2H-2} M^{1/2}(du) M^{1/2}(dv),$$

where $C_H$ is a constant which only depends on $H$.

This conditional covariance can be obtained as the limit $S(s,t) = \lim_{l \to 0^+} S_l(s,t)$, where

$$(2) \qquad S_l(s,t) := C_H \int_0^t \int_0^s |u-v|^{2H-2} l^{-2\psi(1/2)} e^{w_l(u)/2} e^{w_l(v)/2} \, du \, dv.$$

We have the following lemma, which assures this limiting covariance is well defined if (H1/2): $H + \psi(1/2) > 1/2$ holds. Remark that because $\psi(\cdot)$ is convex and $\psi(0) = \psi(1) = 0$, then $\psi(1/2) \leq 0$.

LEMMA 2.2. *Assume* (H1/2) *holds. Then $S(s,t)$ exists a.e. and is not identically equal to zero.*

PROOF. Using the moment generating function of process $\{w_l\}_t$ (Lemma 1 in [7]), for $u \leq v$,

$$\mathbb{E}(e^{1/2(w_l(u)+w_l(v))}) = e^{2\psi(1)\rho_l(0) + (\psi(1)-2\psi(1/2))\rho_l(v-u)},$$

where, by construction, $\psi(1) = 0$. As $\rho_l(0) = \log(1/l) + 1$, we have

$$\mathbb{E}(S_l(s,t)) = C_H \int_0^t \int_0^s |u-v|^{2H-2} e^{-2\psi(1/2)\rho_l(|u-v|)}.$$

On the other hand, $\rho_l(u) = \log(1/u) + \log(u/l) + 1 - u/l$ for $0 < u < l$, and $\rho_l(u) = \log(1/u)$ for $l < u < 1$, so that, for $u > 0$,

$$e^{-2\psi(1/2)\rho_l(u)} \leq e^{-2\psi(1/2)\log(1/u)}$$

for all $l > 0$. An application of the Dominated Convergence Theorem yields that the limiting random covariance function exists a.e.

For the second part of the lemma, by direct calculation of the expectation as above, $\sup_l \mathbb{E}(M_l^{1/2}([0,t])) < \infty$ and $\mathbb{E}(M^{1/2}([0,t])) < \infty$ since $\psi(1/2) > -1/2$ as follows from (H1/2). As in the proof of Theorem 2 (Lemma 3) in [7], this yields the limiting measure $M^{1/2}$ is nondegenerate ($\mathbb{E}(M^{1/2}([0,t])) = t$). For any $0 < a < s \wedge t$ and $0 < \varepsilon < s \wedge t - a$,

$$S(s,t) > \int_0^a \int_{a+\varepsilon}^s |u-v|^{2H-2} M^{1/2}(du) M^{1/2}(dv)$$
$$\geq \varepsilon^{2H-2} M^{1/2}([0,a]) M^{1/2}([a+\varepsilon, s]),$$

where the last expression is not zero a.e. □

Hence, the conditional covariance is well defined, as well as the process $X_t$.

On the other hand, let $B^H(t)$ be a fractional Brownian motion, with $H > 1/2$ and $\mathbb{E}(B^H(1)) = 1$, which is independent of the MRM $M$. It is also possible to define the MFRW-1/2 process directly, as the limit of

$$X_{l,t} := \int_0^t l^{-\psi(1/2)} e^{w_l(u)/2} \, dB_u^H,$$

as $l \to 0^+$. This integral is understood, conditional to $M$, in terms of the Gaussian construction for deterministic integrands. That is, for each $l > 0$, $X_{l,t}|_M$ is a conditionally centered Gaussian process with covariance $S_l(s,t)$.

Let $X'_t$ be the limiting process whenever it exists. The following result shows that both procedures are equivalent.

LEMMA 2.3. *Assume* (H1/2) *holds, then* $X_t \stackrel{\text{law}}{=} X'_t$.

PROOF. That the finite dimensional distributions of $X_{l,t}$ converge to those of $X_t$ follows directly from the definition so that in order to show the stated result it is enough to see that the sequence is tight in the Skorohod metric. For this, we start by observing that

$$\mathbb{E}(|X_{l,t_1} - X_{l,t_2}|^2)$$
$$= \mathbb{E}\left(C_H \int_{t_1}^{t_2} \int_{t_1}^{t_2} |u-v|^{2H-2} l^{-2\psi(1/2)} e^{w_l(u)/2} e^{w_l(v)/2} \, du \, dv\right),$$

as follows by calculating the conditional expectation.

Hence,

(3) $\quad \mathbb{E}(|X_{l,t_1} - X_{l,t_2}|^2)$

(4) $\quad\quad = C_H \int_{t_1}^{t_2} \int_{t_1}^{t_2} |u-v|^{2H-2} e^{-2\psi(1/2)\rho_l(|u-v|)} \, du \, dv$

$$\leq C_H \int_{t_1}^{t_2} |u-v|^{2H-2+2\psi(1/2)} \, du \, dv$$

(5) $$= C(H)|t_1 - t_2|^{2(H-1+\psi(1/2))+2},$$

for a certain constant $C(H)$ and $t_2, t_1 \in [0,1]$, which proves tightness under the stated condition. $\square$

### 2.1.2. A second approach.

DEFINITION 2.2. Let $H > 1/2$. The MFRW-1 is the centered process $X_t$ which, conditionally to $M(dt)$, has a Gaussian distribution with conditional covariance

$$R(s,t) = C_H \int_0^t \int_0^s |u-v|^{2H-2} M(du) M(dv),$$

where $C_H$ is a constant which only depends on $H$.

This conditional covariance can be obtained as the limit $R(s,t) = \lim_{l \to 0^+} R_l(s,t)$, where

(6) $$R_l(s,t) := C_H \int_0^t \int_0^s |u-v|^{2H-2} e^{w_l(u)} e^{w_l(v)} \, du \, dv.$$

The existence of this limiting covariance is assured if the following condition (H1): $H - \psi(2)/2 > 1/2$ holds. Convexity of $\psi(\cdot)$ yields $\psi(2) > 0$.

LEMMA 2.4. *Assume* H1 *holds. Then* $R(s,t)$ *exists a.e. and is not identically equal to zero.*

PROOF. The proof follows exactly as that of Lemma 2.2. Nondegeneracy of $M(dt)$ follows from H1, as $H - \psi(2)/2 > 1/2$ yields $\zeta(2) > 1$. $\square$

Hence, the conditional covariance is well defined, as well as the process $X_t$.

On the other hand, let $B^H(t)$ be a fractional Brownian motion, with $H > 1/2$ and $\mathbb{E}(B^H(1)) = 1$, which is independent of the MRM $M$. It is also possible to define the MFRW-1 process directly, as the limit of

$$X_{l,t} := \int_0^t e^{w_l(u)} \, dB_u^H,$$

as $l \to 0^+$. This integral is understood, conditional to $M$, in terms of the Gaussian construction for deterministic integrands. That is, for each $l > 0$, $X_{l,t}|_M$ is a conditionally centered Gaussian process with covariance

$$R_l(s,t) = C_H \int_0^t \int_0^s |u-v|^{2H-2} e^{w_l(u)} e^{w_l(v)} \, du \, dv.$$

Let $X'_t$ be the limiting process whenever it exists. The following result shows that both procedures are equivalent.

LEMMA 2.5. *Assume* H1 *holds, then* $X_t \stackrel{\text{law}}{=} X'_t$.

PROOF. The proof follows exactly as that of Lemma 2.3. □

2.1.3. *Existence of moments.* An important issue is stating sufficient conditions for the existence of the moments of the MFRW-1/2 or the MFRW-1 process $X_t$. In each case these will be related to the scaling function $\psi(\cdot)$.

Introduce the following definitions:

DEFINITION 2.3. We will say process satisfies assumption H1/2-$p$ if the MFRW-1/2 process $X_t$ satisfies
$$\mathbb{E}(|X_1 - X_0|^p) < \infty.$$

DEFINITION 2.4. We will say process satisfies assumption H1-$p$ if the MFRW-1 process $X_t$ satisfies
$$\mathbb{E}|X_1 - X_0|^p < \infty.$$

Although both constructions are roughly equivalent, the scaling factor will not be the same in every case. Indeed, let $X_t$ be a MFRW-1/2. If $1 \leq j \leq n$, the scaling property (4, in Lemma 2.1) yields
$$\mathbb{E}(|X_{(j-1)/n} - X_{j/n}|^{2p}) = n^{\psi(p) - 4p\psi(1/2) - 2pH} \mathbb{E}(|X_1 - X_0|^{2p}).$$

If, however, $X_t$ is a MFRW-1 process, then we obtain the more familiar
$$\mathbb{E}(|X_{(j-1)/n} - X_{j/n}|^{2p}) = n^{\psi(2p) - 2pH} \mathbb{E}(|X_1 - X_0|^{2p}).$$

In what follows we will restrict our attention to the second approach.

Set $\delta_m := \psi(m) + \psi(m-2) - 2\psi(m-1)$ for $m \geq 2$. Convexity of function $\psi$ assures that $\delta_m \geq 0$ for each $m \geq 2$. Remark that in the case $m = 2$, $\delta_2 = \psi(2)$. The next lemma gives sufficient conditions for H1-2$p$. Its proof is based on power counting techniques, namely, Section 3 (Theorem 3.1) in [12]. Power counting techniques provide sufficient conditions for the existence of integrals whose integrands are products of powers of affine functionals. Based on these results we have the following lemma.

LEMMA 2.6. *Assume* $2H - 1 - \delta_m > 0$ *for each* $2 \leq m \leq 2p$, *then assumption* H1-2$p$ *holds true for* $p \in \mathbb{N}$.

PROOF. Calculating the conditional expectation, we have, for $p \in \mathbb{N}$,

$$\mathbb{E}(|X_{1,1} - X_{l,0}|^{2p}) = c_p \mathbb{E}\left(\int \prod_{i=1}^{p} |u_i - v_i|^{2H-2} e^{w_l(u_i) + w_l(v_i)} du_1 \cdots dv_p\right).$$

Consider any arbitrary permutation of $\{u_1, v_1, \ldots, u_p, v_p\}$ given by $\{z_1, \ldots, z_{2p}\}$. As follows from the properties of the moment generating function of $w_l$, over the set $z_1 < \cdots < z_{2p}$ the last integral is bounded by

$$\int \prod_{i=2}^{p} \prod_{j<i} |z_i - z_j|^{h(i,j)(2H-2) - \delta_{i-j+1}} dz_1 \cdots dz_{2p},$$

where $h(i,j) \in \{0,1\}$ and $h(i,j) = 1$ exactly $p$ times according to whether $z_i = u_k$ and $z_j = v_k$ in the original order [of all the $p(2p-1)$ couples there are exactly $p$ that satisfy this condition]. Following the notation of Section 3 in [12], we are interested in the $p(2p-1)$ linear functionals $L_{i,j} = (z_i - z_j)$ with $j < i$.

By the power counting theorems, existence of the integral will follow if for any strongly independent subset of size $k$ of these linear functionals the condition

$$k + \sum_{l=1}^{k} h(i_l, j_l)(2H - 2) + \delta_{i_l - j_l + 1} > 0$$

holds true. In this setting any strongly independent subset of linear functionals can be identified as a tree (or a collection of disjoint trees) constructed over the set of vertex $\{z_1, \ldots, z_{2p}\}$ identifying each pair $\{z_i, z_j\}$ with an edge. It follows that the maximum size of any strongly independent subset is $2p - 1$.

If $k \leq p$ and if we choose the $k$ terms with $h(i_l, j_l) = 1$, then

$$k + \sum_{l=1}^{k} h(i_l, j_l)(2H - 2) - \delta_{i_l - j_l + 1} = \sum_{l=1}^{k} [(2H - 1) - \delta_{i_l - j_l + 1}] > 0$$

under the assumption that $(2H - 1) - \delta_m > 0$ for each $2 \leq m \leq 2p$.

If $k > p$ or if we choose a subset with $h(i_l, j_l) = 0$, the result follows as $1 - \delta_m > 2H - 1 - \delta_m$. □

REMARK 2.1. When $H = 1/2$, H1-$2p$ is just A$2p$ and is assured whenever $2p - 1 - \psi(2p) > 0$. If $p = 2$, H1-2 is assured if $2H - 1 - \psi(2) > 0$. It is hence fair to conjecture that a sufficient condition for H1-$2p$ is given by $2pH - 1 - \psi(2p) > 0$, which in turn assures together with $2H - 1 - \psi(2) > 0$, that $kH - 1 - \psi(k) > 0$ for any $2 \leq k \leq 2p$, as follows by the properties of function $\psi$.

It is important to remark that in any case $X_t$ is not the subordinated process $U_t = B^H(M[0,t])$. The latter indeed is a centered, conditionally Gaussian process, with conditional covariance

$$\text{Cov}(U_t, U_s) = \tfrac{1}{2}[(M[0,t])^{2H} + (M[0,s])^{2H} - |M[0,t] - M[0,s]|^{2H}].$$

This process is referred to as fractional Brownian motion in multifractal time [13, 16].

**3. Quadratic and higher-order moments of the increments.** Formally, we may be tempted to assume $dX_t = \lim_{l \to 0^+} e^{w_l(t)} dB_t^H$, however, $w_l$ is not necessarily regular and this interpretation may lead to confusion if trying to interpret the above expression by considering pathwise integrals. We will thus consider $X_t$ in terms of the conditionally Gaussian specification. In particular, throughout the rest of the paper we will assume H1 holds true.

Throughout the next sections we will consider dyadic increments, so always $m_n = 2^n$. We also assume, in order to simplify notation, that $T = 1$. Assume we observe $X_{t_1}, \ldots, X_{t_{m_n}}$ with $t_j = j/m_n$. As mentioned in Section 1, usual procedure for estimating $\zeta(p)$ is to look at the average of the $p$ moments of the increments (or, more generally, the increments of order $r$, which are roughly equivalent to considering wavelet coefficients for a wavelet with $r$ vanishing moments) and then use the scaling properties (see [1]). Our main interest in this section will be developing an insight as to the estimation problem at hand, by looking at the asymptotic behavior of these averages for even $p$ and $r = 1$. Because of the long range dependence property of fBm for $H > 1/2$, it is necessary to further restrict the value of the Hurst parameter by asking that $H < 3/4$ in order to obtain a limiting conditional Gaussian behavior. This restriction can be avoided by considering the $r$ increments of $X(t)$, with $r > 1$. However, this involves further technicalities which were not the main interest of this work.

We start by defining the increments of $X_t$. For $j = 1, \ldots, m_n$, define $\Delta X_j = X_{j/m_n} - X_{(j-1)/m_n}$ with $X_0 = 0$.

Set $a_{j,n} = \mathbb{E}_M[(\Delta X_j)^2]^{1/2}$, and define

$$Z_n(p) = \frac{1}{\sqrt{m_n}} \sum_{j=1}^{m_n} [|\Delta X_j|^p - c_p a_{j,n}^p],$$

where $c_p = \mathbb{E}(|Z|^p)$, $Z \sim N(0,1)$.

Finally, define the conditional correlation

$$\rho(k,j) = \rho(\Delta X_j, \Delta X_k) = \frac{\mathbb{E}_M(\Delta X_j \Delta X_k)}{a_{j,n} a_{k,n}}.$$

We are interested in studying the asymptotic behavior of $Z_n(p)$, based on the scaling property and the conditionally Gaussian structure.

By construction,
$$a_{j,n} = \left[\int_{t_{j-1}}^{t_j} \int_{t_{j-1}}^{t_j} |u-v|^{2H-2} M(du) M(dv)\right]^{1/2}$$

and
$$\rho(j,k) = \frac{\int_{t_{j-1}}^{t_j} \int_{t_{k-1}}^{t_k} |u-v|^{2H-2} M(du) M(dv)}{a_{j,n} a_{k,n}}.$$

As mentioned above, the asymptotic behavior of $Z_n$ relies on the conditionally Gaussian behavior of process $X_t$. The next two technical lemmas provide the necessary tools.

LEMMA 3.1. *There exists a (deterministic) constant C, such that*
$$\rho(j,k) \leq C|j-k|^{2H-2}, \tag{7}$$
*for $j \neq k$.*

PROOF. Assume to begin with that $|j-k| \geq 2$. Remark in $\rho(j,k)$, the numerator is bounded by
$$\left(\frac{k-j-1}{m_n}\right)^{2H-2} \left[\int_{t_{j-1}}^{t_j} M(du) \int_{t_{k-1}}^{t_k} M(du)\right].$$

On the other hand, we may bound $a_{l,n}$ from below
$$a_{l,n} > \frac{1}{m_n^{2H-2}} \int_{t_{l-1}}^{t_l} M(du), \tag{8}$$

for $1 \leq l \leq m_n$, which yields the stated bound. For $|j-k|=1$, we have $\rho(j,k) \leq 1$ since $\rho$ is a correlation. □

Define for $r \in \mathbb{N}, r \geq 2$
$$\Gamma_n(r,p) = \frac{1}{m_n} \sum_{1 \leq k, j \leq m_n} \rho(k,j)^r a_{j,n}^p a_{k,n}^p.$$

LEMMA 3.2. *Assume $1/2 < H < 3/4$ and $p \in \mathbb{N}$. Then, $m_n^{2pH-\psi(2p)} \Gamma_n(r,p)$ converges a.s. to a positive integrable r.v. $\Gamma(r,p)$ for all even $r \leq p$, whenever $p$ is even, H1-2p holds and condition A(p): $2p\psi'(2p) < \psi(2p)+1$ is satisfied.*

PROOF. By stationarity, $\mathbb{E}(a_{j,n}^p a_{k,n}^p) \leq \mathbb{E}(a_{0,n}^{2p})$. And as follows by the scaling property, this last expectation is equal to $m_n^{\psi(2p)-2pH} \mathbb{E}(|X_1 - X_0|^{2p})$, which justifies the normalizing constant.

On the other hand, Lemma 3.1 guarantees that $\rho(k,j)^r/C^r|j-k|^{r(2H-2)} \leq 1$.

Set $a_{-j,n} = a_{j,n}$ and define, for $|u| < m_n$,

$$W_{u,n} = \frac{m_n^{2pH-\psi(2p)}}{m_n - |u|} \sum_{k=1}^{m_n-|u|} \frac{\rho(k,k+u)^r}{\max(|u|,1)^{r(2H-2)}} a_{k,n}^p a_{k+u,n}^p.$$

Hence,

$$m_n^{2pH-\psi(2p)} \Gamma_n(r,p) = \sum_{|u| \leq m_n - 1} \left(1 - \frac{|u|}{m_n}\right) [|u|^{r(2H-2)}] W_{u,n}.$$

We will start by showing that $W_{u,n}$ converges a.s. to a certain integrable r.v. $W_{u,r}$ for each fixed $u$, $|u| < m_n$.

I. As part of the proof, we will require a truncation argument which holds for all $|u| < m_n$ simultaneously. Let $r > 1$ be such that $a' := \psi(2pr)/r - \psi(2p) < \sqrt{2 - 1/r} - 1$. This can be seen to hold for $r$ close enough to 1, by checking that $F(1) = 0$ and $F'(1) < 0$ for $F(r) = \psi(2pr) - r\psi(2p) - r\sqrt{2 - 1/r} + r$ under A($p$). Since $\sqrt{2 - 1/r} - 1 < 1 - 1/r$, this choice of $a'$ yields that $a := \psi(2pr)/r - \psi(2p) + 1/r < 1$. Set $\gamma = (1/4 + 3a/4)$, so that $a < \gamma < 1$. Set $A_j = \{m_n^{2pH-\psi(2p)} a_{j,n}^{2p} < m_n^\gamma\}$ and $A_n = \bigcap_{i=1}^{m_n} A_j$.

Then, by the scaling property,

$$P(A_n^c) \leq \sum_{j=1}^{m_n} P(A_j^c) \leq \sum_{j=1}^{m_n} \frac{\mathbb{E}(m_n^{r(2pH-\psi(2p))} a_{j,n}^{2pr})}{n^{r\gamma}}$$

$$\leq m_n^{1+\psi(2pr)-r\psi(2p)-r\gamma} = 2^{nr(a-1)/4},$$

which is the term of a convergent series. Thus, by the Borel–Cantelli lemma, $P(A_n^c \text{ i.o.}) = 0$.

II. To show a.s. convergence of $W_{u,n}$, we start with $u = 0$. Remark in this case $\frac{\rho(k,k)^r}{|\max(0,1)|^{r(2H-2)}} = 1$. The proof follows from a series of steps. The ideas behind this part of the proof are closely related to the proof of Theorem 3.1 in [18].

*Step 1.* Set $l_n = 2^{-(1-\alpha)n}$, where $\alpha = (1-a)/2$, and let $\mathcal{F}_n = \sigma\{w_s(u), s \geq l_n, t \in [0,1]\}$. We shall see $\mathbb{E}_{\mathcal{F}_n}(W_{0,n}) - W_{0,n} \to 0$ a.s.

Set $\Delta_{j,n} = [(j-1)/m_n, j/m_n]$, let $z = (u_1, v_1, \ldots, u_p, v_p)$ and define, for any $0 < l$,

$$d_{n,l}(z) = \frac{m_n^{2pH-\psi(2p)}}{m_n} \sum_j \prod_{i=1}^p |u_i - v_i|^{2H-2} e^{w_l(u_i) - w_{l_n}(u_i)} e^{w_l(v_i) - w_{l_n}(v_i)}$$

$$\times \mathbb{E}\left(\prod_{i=1}^p e^{w_{l_n}(u_i)} e^{w_{l_n}(v_i)}\right) \prod_{i=1}^p \mathbb{1}_{\Delta_{j,n}^2}(u_i, v_i),$$

where $\mathbb{1}_A$ is the indicator function of a set $A \in [0,1]^2$. Also set

$$e_{n,l}(z) = \frac{m_n^{2pH-\psi(2p)}}{m_n} \sum_j \prod_{i=1}^p |u_i - v_i|^{2H-2} \mathbb{1}_{\Delta_{j,n}^2}(u_i, v_i)$$
$$\times \mathbb{E}\left(\prod_{i=1}^p e^{w_l(u_i)} e^{w_l(v_i)}\right).$$

It is important to remark that if $l < l_n$, for each $n$, $w_l(s) - w_{l_n}(s)$ is independent of $w_{l_n}(t)$, for each $s, t \in [0,1]$, because, by construction, the sets are disjoint. Hence,

$$\mathbb{E}\left(\prod_{i=1}^p e^{w_l(u_i)} e^{w_l(v_i)}\right)$$
$$= \mathbb{E}\left(\prod_{i=1}^p e^{w_{l_n}(u_i)} e^{w_{l_n}(v_i)}\right) \mathbb{E}\left(\prod_{i=1}^p e^{w_l(u_i) - w_{l_n}(u_i)} e^{w_l(v_i) - w_{l_n}(v_i)}\right)$$

and $e_{n,l}(z) = \mathbb{E}(d_{n,l}(z))$ so that $\mathbb{E}(\int_A d_{n,l}(z)) = \int_A e_{n,l}(z)$ for any Borel set $A$.

Moreover, let $e_0(p) := \lim_{l \to 0^+} \int e_{n,l}(z) \, dz$, so that

$$e_0(p) = \mathbb{E}\left(\int \prod_{i=1}^p |u_i - v_i|^{2H-2} \, dM(z)\right),$$

as follows by the scaling property and Lemma 2.5 whenever H1 holds.

On the other hand, because the sequence $m_n$ is dyadic, $e^{w_l(u) - w_{l_n}(u)}$ and $e^{w_l(v) - w_{l_n}(v)}$ are independent if $l < l_n$ and $u \in \Delta_{j,n}$ and $v \in \Delta_{k,n}$ with $|j - k| > 2^{n\alpha}$.

Set

$$h_n(z) := \frac{\prod_{i=1}^p e^{w_{l_n}(u_i) + w_{l_n}(v_i)}}{\mathbb{E}(\prod_{i=1}^p e^{w_{l_n}(u_i) + w_{l_n}(v_i)})}.$$

With the above notation, note that $W_{0,n} = \lim_{l \to 0^+} \int h_n(z) d_{n,l}(z) \, dz$ and $\mathbb{E}_{\mathcal{F}_n}(W_{0,n}) = \lim_{l \to 0^+} \int h_n(z) e_{n,l}(z) \, dz$. Now consider

$$\Delta_n := \lim_{l \to 0^+} \int h_n(z)(d_{n,l}(z) - e_{n,l}(z)) \, dz$$

and

$$\tilde{\Delta}_n := \mathbb{1}_{A_n} \lim_{l \to 0^+} \int h_n(z)(d_{n,l}(z) - e_{n,l}(z)) \, dz.$$

First, since $P(\Delta_n \neq \tilde{\Delta}_n) = P(A_n^c)$, it is enough by I to study $\tilde{\Delta}_n$.

Second, by calculating conditional expectations and recalling the terms in the sum are conditionally $2^{n\alpha}$-dependent, we have

$$\mathbb{E}(\tilde{\Delta}_n^2) \leq C \frac{m_n^{\gamma+\alpha}}{m_n^2} \sum_{j=1}^{m_n} m_n^{2pH-\psi(2p)} \mathbb{E}(a_j^{2p}) = O(m_n^{\gamma-1+\alpha}),$$

which is the term of a convergent series in $n$ by the choice of $\gamma$ and $\alpha$.

*Step 2.* We shall now prove that $\lim_{l \to 0^+} \int h_n(z) e_{n,l}(z)\, dz$ converges a.s.

First remark that, for each fixed $l$, $w_l(u) - w_l(0) = P(A_l(u) - A_l(0)) - P(A_l(0) - A_l(u))$, where $A - B = A \cap B^c$ for any two sets $A$ and $B$. By direct calculation, recalling $T = 1$,

$$\int_{A_l(u)-A_l(0)} \mu(dl, ds) = \int_{A_l(0)-A_l(u)} \mu(dl, ds) = u/l + u$$

and

$$\mathbb{E}(e^{q(P(A_l(u)-A_l(0)))}) = \mathbb{E}(e^{q(P(A_l(0)-A_l(u)))}) = e^{\psi(q)(u/l+u)}.$$

Hence, by stationarity, if $\mathbb{E}(|w_{l_n}(0)|^q) < \infty$, there exists $C > 0$ such that, if $u \in \Delta_{j,n}$, then $\mathbb{E}(|w_{l_n}(u) - w_{l_n}(j/m_n)|^q) \leq C(u - j/m_n)/l_n \leq C m_n^{-\alpha}$.

On the other hand, calculating the expectation and the limit as $l \to 0^+$,

$$\lim_{l \to 0^+} \int h_n(z) e_{n,l}(z)\, dz$$

$$= \sum{}' \frac{m_n^{2pH-\psi(2p)}}{m_n}$$

$$\times \sum_{j=1}^{m_n} \int h_n(z) \prod_{i=1}^{p} \mathbb{1}_{\Delta_{j,n}^2}(u_i, v_i)|u_i - v_i|^{2H-2} \prod_{i,k} |z_i - z_k|^{-\alpha(i,k)}\, dz,$$

where $\sum'$ is the sum over all the possible $(2p)!$ orders of the $2p$ integrands. Let $s$ stand for a fixed permutation, such that $s(z_i), s(z_k)$, $i, k = 1, \ldots, p$, are the given positions in the order. Then $\alpha(i,k) = \psi(|s(z_i) - s(z_k)|) + \psi(|s(z_i) - s(z_k)| + 2) - 2\psi(|s(z_i) - s(z_k)| + 1)$ and $\sum_{i,k} \alpha(i,k) = \psi(2p)$ (see [7]).

Consider the expression

$$U_n := \sum{}' \frac{m_n^{2pH-\psi(2p)}}{m_n} \sum_{j=1}^{m_n} e^{2pw_{l_n}(j/m_n)}$$

$$\times \int \frac{\prod_{i=1}^{p} \mathbb{1}_{\Delta_{j,n}^2}(u_i, v_i)|u_i - v_i|^{2H-2} \prod_{i,k} |z_i - z_k|^{-\alpha(i,k)}}{\mathbb{E}(\prod_{i=1}^{p} e^{w_{l_n}(u_i)+w_{l_n}(v_i)})}\, dz,$$

where, for each $j$,

$$\sum{}' m_n^{2pH-\psi(2p)} \int \prod_{i=1}^{p} \mathbb{1}_{\Delta_{j,n}^2}(u_i, v_i)|u_i - v_i|^{2H-2}$$

$$\times \prod_{i,k} |z_i - z_k|^{-\alpha(i,k)} \, dz = e_0(p).$$

And since there exists $c > 0$ such that (Lemma 1 in [7]) $\mathbb{E}(\prod_{i=1}^{p} e^{w_{l_n}(u_i) + w_{l_n}(v_i)}) \geq c m_n^{(1-\alpha)\psi(2p)}$, the expression

$$e_0'(p) := \sum{}' m_n^{2pH - \alpha\psi(2p)}$$

$$\times \int \frac{\prod_{i=1}^{p} \mathbb{1}_{\Delta_{j,n}^2}(u_i, v_i)|u_i - v_i|^{2H-2} \prod_{i,k} |z_i - z_k|^{-\alpha(i,k)}}{\mathbb{E}(\prod_{i=1}^{p} e^{w_{l_n}(u_i) + w_{l_n}(v_i)})} \, dz < \infty$$

and does not depend on $n$, so that we can write

$$U_n = e_0'(p) m_n^{-(1-\alpha)\psi(2p)} \int_0^1 \sum_j e^{2p w_{l_n}(j/m_n)} \mathbb{1}_{\Delta_{j,n}}(u) \, du$$

and by Proposition 2.2 in [18], under A($p$), $U_n$ converges a.s. to a certain random measure $e_0'(p) M^{(2p)}[0, 1]$.

In order to end the proof of this step, set $D_n = \lim_{l \to 0^+} \int h_n(z) e_{n,l}(z) \, dz - U_n$, which can be written as

$$\sum{}' \frac{m_n^{2pH - \psi(2p)}}{m_n} \sum_{j=1}^{m_n} \int \frac{e^{2p w_{l_n}(j/m_n)} - \prod_{i=1}^{p} e^{w_{l_n}(u_i) + w_{l_n}(v_i)}}{\mathbb{E}(\prod_{i=1}^{p} e^{w_{l_n}(u_i) + w_{l_n}(v_i)})}$$

$$\times \prod_{i=1}^{p} \mathbb{1}_{\Delta_{j,n}^2}(u_i, v_i)|u_i - v_i|^{2H-2} \prod_{i,k} |z_i - z_k|^{-\alpha(i,k)}.$$

Set

$$b_{k,n} = \frac{e^{(k-1)w_{l_n}(j/m_n)} \prod_{i=k+1}^{2p} e^{w_{l_n}(z_i)}}{\mathbb{E}(\prod_{i=1}^{p} e^{w_{l_n}(u_i) + w_{l_n}(v_i)})}.$$

Using the identity $z_1 \ldots z_p - w_1 \ldots w_p = (z_1 - w_1)(z_2 \ldots z_m) + w_1(z_2 \ldots z_p - w_2 \ldots w_p)$ and bounding

$$|e^{w_{l_n}(u)} - e^{w_{l_n}(j/m_n)}| \leq e^{w_{l_n}(u)} |1 - e^{w_{l_n}(u) - w_{l_n}(j/m_n)}|$$

$$\leq e^{w_{l_n}(u)} |w_{l_n}(u) - w_{l_n}(j/m_n)|,$$

we have

$$|D_n| \leq \sum{}' \sum_{k=1}^{2p-1} \frac{m_n^{2pH - \psi(2p)}}{m_n}$$

$$\times \sum_{j=1}^{m_n} \int |w_{l_n}(z_k) - w_{l_n}(j/m_n)|$$

$$\times e^{w_{l_n}(j/m_n)} b_{k,n}(z) \prod_{i=1}^{p} \mathbb{1}_{\Delta_{j,n}^2}(u_i, v_i)|u_i - v_i|^{2H-2}$$

$$\times \prod_{i,k} |z_i - z_k|^{-\alpha(i,k)} dz.$$

Let $r > 1$ be as in the proof of step 1. Again, by Lemma 1, [7],

$$\mathbb{E}(e^{rw_{l_n}(j/m_n)} b_{k,n}^r(z)) = O(m_n^{(1-\alpha)(\psi(2pr) - r\psi(2p))}).$$

Let $r' = r/(r-1)$. Hölder's inequality yields

$$m_n^{2pH - \psi(2p)} \int \mathbb{E}(|w_{l_n}(z_k) - w_{l_n}(j/m_n)| e^{w_{l_n}(j/m_n)} b_{k,n}(z))$$

$$\times \prod_{i=1}^{p} \mathbb{1}_{\Delta_{j,n}^2}(u_i, v_i)|u_i - v_i|^{2H-2} \prod_{i,k} |z_i - z_k|^{-\alpha(i,k)} dz$$

$$\leq m_n^{2pH - \psi(2p)} \int [\mathbb{E}(e^{rw_{l_n}(j/m_n)} b_{k,n}^r(z))]^{1/r}$$

$$\times [\mathbb{E}(|w_{l_n}(z_k) - w_{l_n}(j/m_n)|^{r'})]^{1/r'}$$

$$\times \prod_{i=1}^{p} \mathbb{1}_{\Delta_{j,n}^2}(u_i, v_i)|u_i - v_i|^{2H-2}$$

$$\times \prod_{i,k} |z_i - z_k|^{-\alpha(i,k)} dz$$

$$\leq C m_n^{(1-\alpha)(\psi(2pr) - r\psi(2p))/r} m_n^{-\alpha/r'}$$

$$= O(m_n^{(1-\alpha)(\psi(2pr) - r\psi(2p))/r - \alpha(r-1)/r}).$$

The exponent in the last expression can be rewritten as $\eta := (1-\alpha)a' - \alpha(r-1)/r$. In order to end the proof, it is enough to show that $\eta < 0$ and, thus, that the latter is the term of a convergent series.

As follows from the choice of $\alpha$,

$$\eta = a' - \frac{(1 - 1/r - a')^2}{2(a' + 1/r)}.$$

Consider the second degree polynomial $p(x) = 2x(x + 1/r) - (1 - 1/r - x)^2 = (x + \sqrt{2 - 1/r} + 1)(x + 1 - \sqrt{2 - 1/r})$. Hence, $\eta = 2(a' + 1/r)p(a') < 0$ since $0 < a' < \sqrt{2 - 1/r} - 1$.

III. Next, for $u \neq 0$. For this part we require $p$ to be even. We are interested in the convergence of $\frac{m_n^{2pH - \psi(2p)}}{m_n} \sum_k \rho^r(k, k+u) a_k^p a_{k+u}^p$. Since $r$ and $p$ are even, each of the terms in the sum can be written as the product of $2p$

integrals

$$a_k^{p-r} a_{k+u}^{p-r} \left[ \prod_{i=1}^{r} \int_{\Delta_{k,n}} \int_{\Delta_{k+u,n}} |u_i - v_i|^{2H-2} M(du_i) M(dv_i) \right]$$

$$= \prod_{i=1}^{r} \int_{\Delta_{k,n}} \int_{\Delta_{k+u,n}} |u_i - v_i|^{2H-2} M(du_i) M(dv_i)$$

$$\times \prod_{i=1}^{(p-r)/2} \int_{\Delta_{k,n}^2} |u_i - v_i|^{2H-2} M(du_i) M(dv_i)$$

$$\times \prod_{i=(p-r)/2+1}^{(p-r)} \int_{\Delta_{k+u,n}^2} |u_i - v_i|^{2H-2} M(du_i) M(dv_i).$$

This product has exactly the same structure as that of the terms in $W_{0,n}$, although, of course, the expectation will not be the same.

On the other hand, applying the identity, $2ab \leq a^2 + b^2$, to the last product of integrals,

$$\frac{m_n^{2pH-\psi(2p)}}{m_n} \sum_k \rho^r(k, k+u) a_k^p a_{k+u}^p$$

can be seen to be bounded over the set $A_n$ defined in I and the proof of step 1 can be carried out analogously. For step 2, remark that

$$\mathbb{E}(m_n^{2pH-\psi(2p)} \rho^r(k, k+u) a_k^p a_{k+u}^p)$$

does not depend neither on $k$ nor on $n$. The first assertion follows by stationarity. For the second, it is enough to consider a change of variables $z_i' = m_n z_i$. By Lemma 1 in [7], we have, for each fixed $l > 0$,

$$\mathbb{E}\left( \prod_{i=1}^{2p} e^{w_l(u_i)} \right) = e^{\sum_{i=1}^{2p} \sum_{j \leq i} \alpha(i,j) \rho_l(u_i - u_j)},$$

over the set $u_1 \leq \cdots \leq u_{2p}$, where $\rho_l(u) = \log(1/u) + \log(u/l) + 1 - u/l$ for $0 < u < l$ and $\rho_l(u) = \log(1/u)$ for $l < u < 1$. Also, as $\sum_{i,j} \alpha(i,j) = \psi(2p)$, for each fixed $l > 0$, calculating expectations, we find a normalization factor $m_n^{-2pH+\psi(2p)+2p}$, so that by the change of variables, we obtain the factor $m_n^{-2pH+\psi(2p)}$. Hence, by calculating the limit expectation as $l \to 0^+$,

$$e_u(r, p) := \mathbb{E}(m_n^{2pH-\psi(2p)} \rho^r(k, k+u) a_k^p a_{k+u}^p)$$

$$= \mathbb{E} \prod_{i=1}^{r} \int_{[0,1]} \int_{[u,u+1]} |u_i - v_i|^{2H-2} M(du_i) M(dv_i)$$

$$\times \prod_{i=1}^{(p-r)/2} \int_{[0,1]^2} |u_i - v_i|^{2H-2} M(du_i) M(dv_i)$$

$$\times \prod_{i=(p-r)/2+1}^{(p-r)} \int_{[u,u+1]^2} |u_i - v_i|^{2H-2}$$

$$\times M(du_i) M(dv_i).$$

Set

$$e_u(r,p)' := m_n^{2pH - \alpha \psi(2p)} \mathbb{E} \int \prod_{i=1}^{r} \mathbb{1}_{\Delta_{j,n}^2} \mathbb{1}_{\Delta_{j+u,n}^2} |u_i - v_i|^{2H-2}$$

$$\times \prod_{i=r+1}^{r+(p-r)/2} \mathbb{1}_{\Delta_{j,n}^2} |u_i - v_i|^{2H-2}$$

$$\times \prod_{i=r+(p-r)/2+1}^{p} \mathbb{1}_{\Delta_{j+u,n}^2} |u_i - v_i|^{2H-2}$$

$$\times \frac{1}{\mathbb{E}(\prod_{i=1}^{p} e^{w_{l_n}(u_i) + w_{l_n}(v_i)})} M(dz).$$

As in step 2, this expression is bounded and does not depend on $n$. Also, since there are $r$ differences with $|u_i - v_i| > u/m_n$, there exists a constant $C$ such that for all $u$ $\frac{e_u(r,p)'}{|u|^{r(2H-2)}} \leq C$.

Define

$$U_n(u) = \frac{e_u(r,p)' m_n^{-(1-\alpha)\psi(2p)}}{|u|^{r(2H-2)}}$$

$$\times \sum_j \int_0^1 \sum_j e^{2pw_{l_n}(j/m_n)} \mathbb{1}_{\Delta_{j,n}}.$$

Under A($p$), $U_n(u)$ converges a.s. to a certain random measure $W_{u,r} := \frac{e_u(r,p)'}{|u|^{r(2H-2)}} \times M^{(2p)}([0,1])$. Moreover, as in step 2 of II, the random part of this measure does not depend on $u$. Also, define

$$D_n(u) = \mathbb{E}_{\mathcal{F}_n}(W_{0,n}) - U_n(u).$$

In this case bounding the expectation of

$$\mathbb{E}(|w_{l_n}(v) - w_{l_n}(j/m_n)|) \leq C \sum_{k=1}^{u-1} \mathbb{E}(|w_{l_n}((j+k)/m_n) - w_{l_n}(j/m_n)|)$$

$$+ \mathbb{E}(|w_{l_n}(v) - w_{l_n}((j+u-1)/m_n)|)$$

will yield
$$\mathbb{E}(|D_n(u)|) = O(m_n^{(1-\alpha)(\psi(2pr)-r\psi(2p))/r - \alpha(r-1)/r} u).$$

IV. To end the proof, we must show that $m_n^{2pH-\psi(2p)}\Gamma_n(r,p)$ converges a.s. to the limiting r.v. $\Gamma(r,p) = W_0 + 2\sum_{u=1}^{\infty} W_{u,r} u^{r(2H-2)}$. For this, it is enough to consider [as $u^{r(2H-2)}$ is the term of a convergent series]
$$m_n^{2pH-\psi(2p)}\Gamma_n(r,p) - \sum_{u=0}^{m_n} u^{r(2H-2)} U_n(u).$$

Let $A_n$ be as in I. It is enough to show that
$$P\left(\left\{\sum_{u=0}^{m_n} u^{r(2H-2)}|W_{u,n} - U_n(u)| > \varepsilon\right\} \cap A_n\right)$$

is the term of a convergent series. This follows from remarking that, for any $m$, the above probability is bounded by
$$Cm_n^{(\gamma-1+\alpha)/2 + (1-\alpha)(\psi(2pr)-r\psi(2p))/r - \alpha(r-1)/r} m^2$$
$$+ \sup_{u \leq m_n} \mathbb{E}(W_{u,r} + W_{u,n}) m^{r(2H-2)+1}$$

and that $\mathbb{E}(W_{u,r}) \leq C\mathbb{E}(W_{0,r})$ and $\mathbb{E}(W_{u,n}) \leq C\mathbb{E}(W_{0,n})$.

The limit $\Gamma(r,p)$ is integrable because $W_0$ is and is a.s. positive as $\Gamma(r,p) > W_0$. $\square$

Let $g_r$ be the coefficients of the expansion of $G(x) = |x|^p - c_p$ over the Hermite polynomials (since $G$ is a centered even function, $g_r = 0$ for $r = 0, 1$). Set $\Gamma(p) = \sum_r g_r^2 r! \Gamma(r,p)$. Remark that since $\sum_r g_r^2 r! < \infty$ and $\Gamma_n(r,p) \leq \Gamma_n(2,p)$, we have $\Gamma_n(p) := \sum_r g_r^2 r! \Gamma_n(r,p)$ converges a.s. to $\Gamma(p)$, where the latter is a positive integrable r.v.

Lemmas 3.1 and 3.2 assure the conditional convergence of $m_n^{pH-\psi(2p)/2} Z_n(p)$ toward a normal random variable. More precisely, we have the following:

THEOREM 3.1. *If $1/2 < H < 3/4$, then for even $p$, under H1-$2p$ and condition $A(p)$, given in Lemma 3.2, then a.s.*
$$\mathbb{E}_M(e^{i\gamma m_n^{pH-\psi(2p)/2} Z_n(p)}) \to e^{-\gamma^2 \Gamma(p)/2}.$$

PROOF. Write
$$m_n^{pH-\psi(2p)/2} Z_n(p) = \frac{m_n^{pH-\psi(2p)/2}}{\sqrt{m_n}} \sum_{1 \leq k \leq m_n} a_{k,n}^p G(\Delta X_k/a_{k,n}),$$

where, by Lemma 3.1, $\{\Delta X_k/a_{k,n}\}_{k \in \mathbb{N}}$ is a (conditionally) Gaussian and centered sequence of strongly dependent r.v., that is, such that $\text{Cov}(\Delta X_k/a_{k,n},$ $\Delta X_j/a_{j,n}) = O(|k - j|^{2H-2})$. The result follows directly since $G$ has Hermite rank equal to 2, [3, 9] by calculating the conditional moments of $Z_n(p)$ and by the a.s. convergence of $\Gamma_n(r,p)$. Indeed, using moment calculating techniques, it is possible to show that

$$\mathbb{E}_M([m_n^{pH-\psi(2p)/2} Z_n(p)]^{2q}) = c_{2q} \Gamma_n(p)^q (1 + o_n(1))$$

and

$$\mathbb{E}_M([m_n^{pH-\psi(2p)/2} Z_n(p)]^{2q+1}) \leq o_n(1) \Gamma_n(p)^{(2q+1)/2}.$$

The proof then follows by the a.s. convergence of $\Gamma_n(p)$. □

**4. Statistical applications.** In this section we study a moment based estimator of $\zeta(p)$. Moment estimators have the advantage of being very easy to implement. However, experimental studies show the existence of a bias problem as discussed in Section 1.

Based on the results of Section 3, we will discuss rates of convergence for these estimators and offer insights as to some of the estimation problems at hand. Namely, we are interested in the construction of asymptotic confidence intervals. When $m_n \gg T$ ($T = 1$ in order to simplify), we show the existence of a random $1/\log m_n$ bias term, whereas the conditional CLT yields an estimation error of order $1/(m_n^{1/2+\psi(p)-\psi(2p)/2} \log(m_n))$.

When $T \sim m_n$, the scaling property of the MRM yields the appearance of a linear bias term, as shall be discussed in detail in Section 4.2. In this case approximate confidence intervals are given.

Based on the notation introduced in Section 3, $\mathcal{D}_M$ will stand for convergence in distribution conditional to the MRM $M$.

4.1. *Case $T = 1$.* Throughout this section, we will assume $p$ is an even positive integer. Set to begin with

$$S_n(p) = \frac{m_n^{pH-\psi(p)}}{m_n} \sum_{j=1}^{m_n} |\Delta X_j|^p$$

and

$$B_n(p) = \frac{m_n^{pH-\psi(p)}}{m_n} \sum_{j=1}^{m_n} c_p a_{j,n}^p.$$

We have a first result concerning the convergence of $B_n(p)$.

LEMMA 4.1. *Assume* H1-$p$ *holds. Then, $B_n(p)$ converges a.s. to a certain a.s. positive random variable $B(p)$.*

PROOF. Convergence of $B_n(p)$ to the random measure $B(p) = e'_0(p/2) \times M^{(p)}[0,1]$ follows as in the proof of Lemma 3.2. □

Using a first-order Taylor expansion and Theorem 3.1 yield, whenever $\psi(2p) < 2\psi(p) + 1$

$$m_n^{1/2+\psi(p)-\psi(2p)/2} \log\left(\frac{S_n(p)}{B_n}\right) = m_n^{1/2+\psi(p)-\psi(2p)/2} \frac{S_n - B_n(p)}{B(p) + o_n(1)}$$

$$\xrightarrow{\mathcal{D}_M} \frac{N(0,\Gamma(p))}{B(p)}.$$

The above discussion yields the following result.

COROLLARY 4.1. *Assume* H1-2p *holds true for p an even integer. Under the conditions of Theorem 3.1,*

$$\frac{\log S_n(p)}{\log m_n} = b_n + e_n,$$

*where* $\log(m_n)b_n \to \log(B(p))$ *a.s. and* $m_n^{1/2+\psi(p)-\psi(2p)/2} \log(m_n) \frac{B(p)}{\Gamma(p)^{1/2}} e_n \xrightarrow{\mathcal{D}_M} N(0,1)$.

4.2. *Case $T \sim m_n$: a certain linearization effect.* In this subsection we assume that $T = O(m_n)$. That is, that the range of scale invariance of the MRM $M$ is large enough to be comparable to the (fixed) number of observations. In this case the scaling property accounts for a certain linearization phenomena in the estimation procedure.

We are interested in nonasymptotic bounds for $|\log(S_n(p))/\log(m_n) - \mathbb{E}(\log(S_n(p))/\log(m_n))|$. More precisely, we shall see in this case that the scaling property yields $\mathbb{E}(|\log(S_n(p))/\log(m_n) - (H + \psi'(0))p|) \le C/\sqrt{\log m_n}$, assuming $\psi$ is two times differentiable at zero.

We must first introduce some additional notation. Define

$$\tilde{a}_j = \left[\int_{j-1}^{j}\int_{j-1}^{j} |u-v|^{2H-2} M(du)M(dv)\right]^{1/2}$$

and

$$\tilde{\rho}(j,k) = \frac{\int_{j-1}^{j}\int_{k-1}^{k} |u-v|^{2H-2} M(du)M(dv)}{\tilde{a}_j \tilde{a}_k}.$$

By the scaling property, 4 in Lemma 2.1, we have

$$\rho(j,k) \stackrel{\text{law}}{=} \tilde{\rho}(k,j),$$

$$m_n^{2H} a_{k,n} a_{j,n} \stackrel{\text{law}}{=} e^{2\Omega_{1/m_n}} \tilde{a}_k \tilde{a}_j.$$

Next introduce the random matrix $\tilde{\Sigma}_n := [\tilde{a}_k \tilde{a}_j \tilde{\rho}(j,k)]_{j,k}$ and consider $\Sigma_n$ the conditional covariance matrix of $m_n^H(\Delta X_1, \ldots, \Delta X_n)$. The scaling property thus yields

$$\Sigma_n \stackrel{\text{law}}{=} e^{2\Omega_{1/m_n}} \tilde{\Sigma}_n.$$

For each $n$, let $(\Delta \tilde{X}_1, \ldots, \Delta \tilde{X}_n)$ be a conditionally Gaussian centered random vector with conditional covariance matrix $\tilde{\Sigma}_n$. By the scaling property

$$m_n^H(\Delta X_1, \ldots, \Delta X_n) \stackrel{\text{law}}{=} e^{\Omega_{1/m_n}}(\Delta \tilde{X}_1, \ldots, \Delta \tilde{X}_n).$$

Remark that, as follows from the scaling property, the same bounds hold a.s. for $\tilde{\rho}$ as for $\rho$. Define for $r \in \mathbb{N}, r \geq 2$

$$\tilde{\Gamma}_n(r,p) = \frac{1}{m_n} \sum_{1 \leq k,j \leq m_n} \tilde{\rho}(k,j)^r \tilde{a}_j^p \tilde{a}_k^p.$$

Then, by stationarity of the MRM, with the notation introduced in the proof of Lemma 3.2,

$$\mathbb{E}(m_n^{2pH} \tilde{\Gamma}_n(r,p)) \leq C^r e_0(p)\left(1 + 2\sum_{u \geq 1} u^{r(2H-2)}\right),$$

whenever H1-$2p$ holds.

Introduce, for $p \in \mathbb{N}$,

$$\tilde{Z}_n(p) = \frac{1}{\sqrt{m_n}} \sum_{j=1}^{m_n} [|\Delta \tilde{X}_j|^p - c_p \tilde{a}_j^p].$$

By the scaling property,

$$Z_n(p) \stackrel{\text{law}}{=} e^{p\Omega_{1/m_n}} \tilde{Z}_n(p),$$

where, by construction, the random variable $\Omega_{1/m_n}$ is independent both of the collection of random variables $\{M[k,k+1]\}_{k \in \mathbb{Z}}$ for $k \leq m_n$ and the fractional Brownian motion $B_t^H$.

Calculating the conditional variance,

$$\mathbb{E}_M((\tilde{Z}_n(p))^2) = \sum_r g_r^2 r! \tilde{\Gamma}_n(r,p),$$

so that under H1-$2p$, we have

$$\mathbb{E}(m_n^{2pH}(\tilde{Z}_n(p))^2) \leq e_0(p) \sum_r g_r^2 r! C^r \left(1 + 2\sum_{u \geq 1} u^{r(2H-2)}\right).$$

Next consider, for even $p$,

$$\tilde{B}_n(p) = \frac{1}{m_n} \sum_{j=1}^{m_n} c_p \tilde{a}_j^p.$$

For each fixed $n$, $\tilde{B}_n(p)$ is a bounded and positive r.v. under H1-$p$, such that there exists a constant $k_1$, $\mathbb{E}(\tilde{B}_n(p)) \leq k_1 e_0(p/2)$. Jenssen's inequality and positivity of $e_0(p/2)$ thus yield that there exists a constant $k_2$ with $|\mathbb{E}(\log(\tilde{B}_n(p)))| \leq k_2$. Let

$$\tilde{S}_n(p) = \frac{m_n^{pH}}{m_n} \sum_{j=1}^{m_n} |\Delta \tilde{X}_j|^p.$$

Once again the scaling property assures

$$\log\left(\frac{S_n}{B_n}\right) \stackrel{\text{law}}{=} \log\left(\frac{\tilde{S}_n}{\tilde{B}_n}\right).$$

On the other hand,

$$\log(B_n) \stackrel{\text{law}}{=} \log(\tilde{B}_n) + p\Omega_{1/m_n},$$

where $\Omega_{1/m_n}$ is an infinitely divisible random measure, independent of $\tilde{B}_n$ and of the fractional Brownian process $B_t^H$, and such that $\mathbb{E}(e^{q\Omega_\lambda}) = e^{-\psi(q)\log(\lambda)}$. Note that the latter yields

$$\Omega_{\lambda_1 \lambda_2} \stackrel{\text{law}}{=} \Omega_{\lambda_1} + \Omega_{\lambda_2},$$

with $\Omega_{\lambda_1}$ and $\Omega_{\lambda_2}$ independent r.v. Hence, if $m = \log m_n$, writing

$$1/m_n = ((1/m_n)^{1/m})^{[m]}((1/m_n)^{1/m})^{m-[m]},$$

we have

$$\Omega_{1/m_n} \stackrel{\text{law}}{=} \Omega_1' + \cdots + \Omega_{[\log m_n]}' + \Omega_n^\star,$$

where $\Omega_j'$ are i.i.d. r.v. with $\mathbb{E}(e^{q\Omega_j'}) = e^{\psi(q)}$. Let $q$ be such that $\zeta(q) > 1$. By definition, $\mathbb{E}(e^{q\Omega^\star}) = e^{-\psi(q)(m-[m])\log(1/e)}$. By Markov's inequality,

$$P(\Omega^\star > m) \leq \frac{e^{\psi(q)m}}{e^{qm}} = e^{-\zeta(q)m},$$

which is the term of a convergent series and $\Omega_n^\star/\log m_n \to 0$ a.s. It follows ([11], page 565) that $\Omega_{1/m_n}/\log(m_n) \to \psi'(0)$ in probability whenever $\psi'(0)$ exists. Moreover, if $\psi''(0)$ exists, then $\mathbb{E}((\Omega')^2) < \infty$ ([11], page 512), which by the SLLN assures $\Omega_{1/m_n}/\log(m_n) \to \psi'(0)$ a.s. and that

(9) $$\sqrt{\log m_n}\left[\frac{\Omega_{1/m_n}}{\log m_n} - \psi'(0)\right] \to N(0, \sigma_\Omega^2)$$

by the CLT, where $\sigma_\Omega^2 = \text{Var}(\Omega_1')$.

An application of Markov's inequality and the above discussion yield the following result.

COROLLARY 4.2. *Assume $\psi''(0)$ exists and set $\mu := \psi'(0)$. Assume* H1-$2p$ *holds true. Then if $m_n = O(T)$, there exists a positive constant $C$ such that*

$$P\left(\left|\frac{\log S_n(p)}{\log m_n} - (H+\mu)p\right| > \varepsilon\right) \leq \frac{C}{\sqrt{\log(m_n)}\varepsilon}.$$

DEPARTAMENTO DE MATEMÁTICAS
IVIC CARRETERA PANAMERICANA
KM 11, APTDO 21827
CARACAS 1020A
VENEZUELA
E-MAIL: cludena@ivic.ve